\documentclass[12pt,a4paper,reqno]{article}
\usepackage{graphics}
\usepackage{epsfig}

\usepackage{amsmath,amsthm,amssymb}
\newtheorem{theorem}{Theorem}[section]
\newtheorem{corollary}[theorem]{Corollary}

\newtheorem{example}[theorem]{Example}
\theoremstyle{definition}

\theoremstyle{remark}
\newtheorem{remark}[theorem]{Remark}
\numberwithin{equation}{section}

\begin{document}
\title{  Minimum Roman Dominating Distance Energy of a Graph}

\author{Lakshmanan. R \\
Department of Mathematics\\
Central University of Tamil Nadu\\
Tiruvarur-610005, Tamil Nadu, India\\
{Email:laksh1983ram@gmail.com }
\bigskip\\
N. Annamalai\\
Central University of Tamil Nadu\\
Tiruvarur-610005, Tamil Nadu, India\\
{Email: algebra.annamalai@gmail.com}
}
\date{}
\maketitle

\vspace*{0.5cm}
\begin{abstract} 
In this correspondence, we introduced the concept of minimum roman dominating  distance  energy $E_{RDd}(G)$ of a graph $G$ and computed minimum roman dominating  distance  energy of some standard graphs. Also, we discussed the properties of eigenvalues of a minimum roman dominating distance  matrix $A_{RDd}(G).$ Finally, we derived the  upper and lower bounds for $E_{RDd}(G).$ 

\end{abstract}
\vspace*{0.5cm}

{\it Keywords:}   Energy of a graph, Roman dominating function, Roman domination.

{\it 2010 Mathematical Subject Classification:} Primary 05C50, 05C69
\vspace{0.5cm}
\vspace{1.5cm}
\section{Introduction}
In 1978, I. Gutman was introduced the concept of energy of a graph \cite{gut}. The graph $G=(V,E)$ mean a simple connected graph with $n$ vertices and $m$ edges.  The distance between two vertices $u$ and $v$ is the length  shortest distance between $u$ and $v.$
 The Wiener index $W(G)$ of $G,$ is the sum of the lengths of the shortest paths between all pairs of vertices.
Let $A=(a_{ij})$ be the adjacency matrix of a graph $G.$ Then the energy $E(G)$ of a graph $G$ is defined by the sum of absolute value of all eigenvalues of $A.$ For more details about energy of a graph\cite{bapat}.

The union of two simple graphs $G_1 = (V_1, E_1)$ and $G_2= (V_2, E_2)$ is the simple graph with vertex set $V_1 \cup V_2$ and edge set
$E_1\cup E_2.$	The union of $G_1$ and $G_2$ is denoted by $G_1 \cup G_2.$
A {\it crown graph } $S_k^0$ is a bipartite graph with two sets of vertices $\{a_1, a_2, \dots, a_k\}$ and $\{b_1, b_2,\dots, b_k\}$ and with an edge from $a_i$ to $b_j$ whenever $i \neq j.$
For a positive integer $n\geq2,$ a healthy spider is a star $K^*_{1, n-1}$ with all of its edges subdivided\cite{EJcock}.

The distance matrix $A_{d}=(d_{ij})$ of $G$ is a symmetric matrix of order $n$ where $d_{ij}$ is the distance between  $i^{th}$ and $j^{th}$ vertices of a graph.  The distance energy $E_{d}(G)$ of the graph $G$ is defined by the sum of absolute value of all eigenvalues of $A_{d}.$ The distance matrix of an
undirected graph has been widely studied in the literature, see \cite{kir,gra,graham,poll}.

A set $S \subseteq V$ is a dominating set if every vertex of $V\setminus S$
is adjacent to at least one vertex in $S.$ The domination number $\gamma(G)$ is the minimum cardinality of a dominating set in $G,$ and a dominating set $S$ of minimum cardinality is called a $\gamma$-set of $G.$ E.J. Cockayne et al.(2004) introduce the concept of roman domination in graphs\cite{EJcock}.
A {\it Roman dominating function} on a graph $G = (V, E)$ is a function
$f : V \rightarrow \{0, 1, 2\}$ satisfying the condition that every vertex $u$ for which $f(u) = 0$ is adjacent to at least one vertex $v$ for which $f(v) = 2.$

Let $(V_0, V_1, V_2)$ be the ordered partition of $V$ induced by $f,$ where $V_i = \{v \in V\mid f(v) = i\}$ and $| V_i| = n_i,$ for $i =0, 1, 2.$ Note that there exists a one-one correspondence between the functions $f : V \rightarrow \{0, 1, 2\}$ and the ordered partitions $(V_0, V_1, V_2)$ of $V.$ Thus, we will write $f = (V_0, V_1, V_2).$ A function $f = (V_0, V_1, V_2)$ is a Roman dominating function (RDF) if $V_2\succ V_0,$ where $\succ$ means that the set $V_2$ dominates the set $V_0.$ The weight of $f$ is $f(V) =\sum\limits_{v\in V} f(v)=2n_2 + n_1.$

The Roman domination number, denoted $\gamma_R(G),$ equals the minimum weight of an RDF of $G,$ and we say that a function $f = (V_0, V_1, V_2)$ is a $\gamma_R$-function if it is an RDF and $f(V) = \gamma_R(G).$  

\begin{theorem}\cite{EJcock}\label{a}
For any graph $G,$  $\gamma(G) \leq \gamma_R(G)\leq 2 \gamma(G).$
\end{theorem}

M.R. Rajesh kannan et al.(2013) studied the minimum covering distance energy of a graph\cite{prad} and also they were studied Laplacian minimum dominating energy of a graph\cite{dhar}.
 M.R. Rajesh kannan et al.(2014) introduced the concept of the  minimum dominating distance energy of a graph\cite{MRRaj}. Let $G$ be the graph with vertex set $V=\{v_1,v_2, \dots, v_n\}$. Let $D$ be a minimum dominating set of a graph $G$. The minimum dominating distance matrix of $G$ is the square matrix $A_{Dd}(G):=(d_{ij})$ where 
  $$d_{ij}=
 \begin{cases}
1 & \text{ if } i=j \text{ and } v_{i}\in D\\
d(v_i,v_j) & \text{ otherwise }
 \end{cases}$$  
Let $\delta_1, \delta_2, \dots, \delta_n$ be the eigenvalues of $A_{Dd}(G).$ Then the minimum dominating energy $E_{Dd}(G)$ of $G$ is $$E_{Dd}(G)=\sum\limits_{j=1}^{n} |\delta_j|.$$

In this paper, we introduce the concept of minimum roman dominating distance energy of a graph in section 2. In Section 3, we find the minimum roman dominating distance energy of some standard graphs. In section 4, we discussed the properties of eigenvalues of a minimum roman dominating distance matrix $A_{RDd}(G).$ We derived the  upper and lower bounds for $E_{RDd}(G)$ in section 5. 

\section{The Minimum Roman Dominating Energy}
In this section, we introduce the concept of minimum roman dominating energy of a graph.

Let $f=(V_0, V_1, V_2)$ be a $\gamma_R$-function of a graph $G.$ The minimum roman dominating distance matrix $A_{RDd}(G)$ of $G$ is defined as  
$A_{RDd}(G):=(\bar{d}_{ij})$ where 
  $$\bar{d}_{ij}=
 \begin{cases}
2 & \text{ if } i=j \text{ and } v_{i}\in V_2\\
1 & \text{ if } i=j \text{ and } v_{i}\in V_1\\
d(v_i, v_j) & \text{ otherwise }
 \end{cases}$$  
Let $\rho_1, \rho_2, \dots, \rho_n$ be the eigenvalues of $A_{RDd}(G).$ Then the minimum roman dominating distance energy $E_{RDd}(G)$ of $G$ is defined as 
$$E_{RDd}(G)=\sum\limits_{k=1}^{n} |\rho_k|.$$ 
Note that $tr(A_{RDd}(G))=\gamma_R(G).$ 
\begin{example}
The minimum roman dominating function of the following graph $G$ 
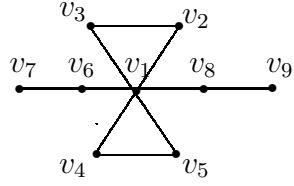
\begin{center}
\begin{picture}(102.499,69)(0,0)
\multiput(33.749,28.75)(.03125,-.03125){8}{\line(0,-1){.03125}}
\put(4.249,42){\line(1,0){95.75}}
\multiput(30.999,65.75)(.03372210953,-.04969574037){986}{\line(0,-1){.04969574037}}
\put(31.249,65.5){\line(1,0){34}}
\put(63.999,17){\line(-1,0){30.25}}
\multiput(64.749,65.5)(-.03370488323,-.05148619958){942}{\line(0,-1){.05148619958}}
\put(4.957,41.958){\circle*{2.915}}
\put(28.457,41.958){\circle*{2.915}}
\put(48.707,40.958){\circle*{2.915}}
\put(73.957,41.958){\circle*{2.915}}
\put(99.707,42.208){\circle*{2.915}}
\put(31.707,65.458){\circle*{2.915}}
\put(64.707,65.458){\circle*{2.915}}
\put(33.957,17.458){\circle*{2.915}}
\put(63.707,17.208){\circle*{2.915}}
\put(6.499,50){\makebox(0,0)[cc]{$v_7$}}
\put(49.499,50){\makebox(0,0)[cc]{$v_1$}}
\put(25,71){\makebox(0,0)[cc]{$v_3$}}
\put(71,68.5){\makebox(0,0)[cc]{$v_2$}}
\put(25,12){\makebox(0,0)[cc]{$v_4$}}
\put(71,12){\makebox(0,0)[cc]{$v_5$}}
\put(73.749,50){\makebox(0,0)[cc]{$v_8$}}
\put(102.499,50){\makebox(0,0)[cc]{$v_9$}}
\put(28.249,50){\makebox(0,0)[cc]{$v_6$}}
\end{picture}
\end{center}
is
$f=(V_0, V_1, V_2)$ where $V_2=\{v_1\}, V_1=\{v_7, v_9\}$ and $V_0=\{v_2, v_3, v_4, v_5, v_6, v_8\}.$
Then the minimum roman dominating distance matrix is
$$A_{RDd}(G)=\left(\begin{array}{rrrrrrrrr}
2 & 1 & 1 & 1 & 1 & 1 & 2 & 1 & 2 \\
1 & 0 & 1 & 2 & 2 & 2 & 3 & 2 & 3 \\
1 & 1 & 0 & 2 & 2 & 2 & 3 & 2 & 3 \\
1 & 2 & 2 & 0 & 1 & 2 & 3 & 2 & 3 \\
1 & 2 & 2 & 1 & 0 & 2 & 3 & 2 & 3 \\
1 & 2 & 2 & 2 & 2 & 0 & 1 & 2 & 3 \\
2 & 3 & 3 & 3 & 3 & 1 & 1 & 3 & 4 \\
1 & 2 & 2 & 2 & 2 & 2 & 3 & 0 & 1 \\
2 & 3 & 3 & 3 & 3 & 3 & 4 & 1 & 1
\end{array}\right)_{9\times 9}$$
Then the characteristic equation of $A_{RDd}(G)$ is $$\rho^{9} - 4\rho^{8} - 171\rho^{7} - 1034\rho^{6} - 2339\rho^{5} - 1284\rho^{4} + 2659\rho^{3} + 4438\rho^{2} + 2410\rho + 444=0$$ and the eigenvalues are $\rho_1=-3, \rho_2=-1, \rho_3=-1, \rho_4\approx-4.5615, \rho_5\approx -0.4384, \rho_6\approx -3.9721, \rho_7\approx-0.8397, \rho_8\approx1.2642,$ and $\rho_9\approx 17.5476.$ Hence the minimum roman dominating energy of $G$ is $E_{RDd}(G)\approx 33.6237.$ 

Note that this graph has unique minimum roman dominating function.
\end{example}

\section{Minimum Roman Dominating Distance Energy of  Some Standard Graphs}
In this section, we studied the minimum roman dominating distance energy of complete, complete bipartite, crown, star  and healthy spider graphs.

Denote $J_n$ is an $n\times n$ all ones matrix, $I_n$ is an $n\times n$ identity matrix, $D_k$ is a diagonal matrix whose $k^{th}$ diagonal entry is zero and other diagonal entries are two, $\bf{a}_{n}$ is a $1\times n$ row vector $[a,a,\dots,a]$ and $\bf{a}'_{n}$ is the transpose of $\bf{a}_{n}.$
\begin{theorem}
For any integer $n\geq 3,$ $E_{RDd}(K_n)=2n-2.$
\end{theorem}
\begin{proof}
For a complete graph $K_n,$ the minimum roman dominating function is $f=(V_0, V_1, V_2)$ where $V_2=\{v_i\}$ for any $i,$ $V_1=\emptyset$ and $V_0=V\setminus V_2.$ Then the minimum roman dominating distance matrix $A_{RDd}(K_n)=(a_{ij})$ where
$$a_{ij}=
 \begin{cases}
2 & \text{ if } i=j \text{ and } v_{i}\in V_2\\
0 & \text{ if } i=j \text{ and } v_{i}\in V_0\\
1 & \text{ otherwise }
 \end{cases}$$
One can easily show that the characteristic polynomial of $A_{RDd}(K_n)$ is  $(\rho+1)^{n-2}(\rho^2-n\rho +n-3).$ Hence the eigenvalues are -1 with multiplicity $n-2$ and $\frac{n\pm \sqrt{n^2-4n+12}}{2}.$ Since for $n\geq 3,$ $n\geq \sqrt{n^2-4n+12}.$ Therefore, the eigenvalues $\frac{n+ \sqrt{n^2-4n+12}}{2}$ and $\frac{n- \sqrt{n^2-4n+12}}{2}$ are positive. Hence the sum of absolute values of all eigenvalues is $2n-2.$ That is, $E_{RDd}(K_n)=2n-2.$
\end{proof}

\begin{corollary}
For any integer $n\geq 3,$ $E_{RDd}(K_n)=E_d(K_n)=E(K_n).$
\end{corollary}
\begin{proof}
Let $n\geq 2$ be an integer. The adjacency matrix of a complete graph $K_n$ is $J_n-I_n$ and the eigenvalues are $n-1$ and $-1$ with multiplicity $n-1,$  the energy of $K_n$ is $2n-2.$
\end{proof}

\begin{theorem}
For any integer $r\geq 2,$ $$E_{RDd}(K_{r,r})=2(2r-4)+\sqrt{(r-2)^2+8}+\sqrt{(3r-2)^2+24}.$$
\end{theorem}
\begin{proof}
Let $X=\{v_1, v_2, \dots, v_r\}$ and $Y=\{w_1, w_2, \dots, w_r\}$ be a partition of the vertex set of a  complete bipartite graph $K_{r,r}.$ Then the minimum roman dominating function is $f=(V_0, V_1, V_2)$ where $V_2=\{v_i, w_j\}$ for any $1\leq i, j\leq r,$ $V_1=\emptyset$ and $V_0=V\setminus V_2.$ Then the minimum roman dominating distance matrix is
$$A_{RDd}(K_{r,r})=\left[
\begin{array}{c|c}
2J_r-D_i & J_r \\ \hline
J_r & 2J_r-D_j
\end{array}\right]_{2r\times 2r}$$
One can easily show that the characteristic equation of $A_{RDd}(K_{r,r})$ is  $(\rho+2)^{2r-4}(\rho^2-(r-2)\rho -2)(\rho^2-(3r-2)\rho-6)=0.$  Then the eigenvalues are $-2$ with multiplicity $2r-4,$ $\frac{(r-2)\pm \sqrt{(r-2)^2+8}}{2}$ and $\frac{(3r-2)\pm \sqrt{(3r-2)^2+24}}{2}.$ Therefore, the eigenvalues $\frac{(r-2)- \sqrt{(r-2)^2+8}}{2},$  $\frac{(3r-2)- \sqrt{(3r-2)^2+24}}{2}$ are negative  and the eigenvalues $\frac{(r-2)+ \sqrt{(r-2)^2+8}}{2},$ $\frac{(3r-2) +\sqrt{(3r-2)^2+24}}{2}$ are positive.

Hence the sum of absolute values of all eigenvalues is $2(2r-4)+\sqrt{(r-2)^2+8}+\sqrt{(3r-2)^2+24}.$ That is, $E_{RDd}(K_{r,r})=2(2r-4)+\sqrt{(r-2)^2+8}+\sqrt{(3r-2)^2+24}.$
\end{proof}

\begin{theorem}
For any $n\geq 3,$ $E_{RDd}(K_{1,n-1})=4n-6.$
\end{theorem}
\begin{proof}
Consider the star graph $K_{1, n-1}$ with vertex set $V=\{v_0, v_1, v_2, \dots, v_{n-1}\},$ where $deg(v_0)=n-1.$ The minimum roman dominating distance function is $f=(V_0, V_1, V_2)$ where $V_2=\{v_0\},$ $V_1=\emptyset$ and $V_0=V\setminus V_2=\{v_1, v_2,\dots, v_{n-1}\}.$ Then the minimum roman dominating distance matrix is
  $$
  A_{RDd}(K_{1,n-1}) =
  \left[
\begin{array}{c|c}
2&\bf{1}_{n-1}\\ \hline
\bf{1}'_{n-1}&2J_{n-1}-2I_{n-1}  \\ \end{array}\right]_{n\times n}$$
  
The characteristic equation is $(\rho+2)^{n-2}(\rho^2-(2n-2)\rho+3n-7)=0.$
Then the eigenvalues are $-2$ with multiplicity $n-2,$ $(n-1)\pm \sqrt{n^2-5n+8}.$ 
Hence the sum of absolute values of all eigenvalues is $4n-6.$ That is, $E_{RDd}(K_{1,n-1})=4n-6.$
\end{proof}

\begin{theorem}\label{b}
For an odd integer $n\geq 4,$ $$11n-19\leq E_{RDd}(K^{*}_{1,n-1})\leq 6n^2-4n-16.$$
\end{theorem}
\begin{proof}
Consider the healthy spider graph $K^{*}_{1, n-1}$ with vertex set
$V=\{v_0, v_1, v_2, \dots, v_{n-1},\\ u_1, u_2, \dots, u_{n-1}\}.$ The vertex $v_0$ is adjacent with $v_1, v_2, \dots, v_{n-1}$ and for $1\leq i\leq n-1,$ $u_i$ is adjacent with $v_i.$ Then the minimum roman dominating  function is $f=(V_0, V_1, V_2)$ where $V_2=\{v_0\},$ $V_1=\{u_1, u_2, \dots,u_{n-1}\}$ and $V_0=\{v_1, v_2,\dots, v_{n-1}\}.$ 
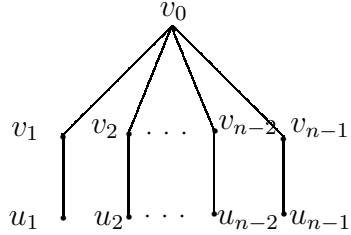
\begin{center}
\begin{picture}(105.75,82)(0,0)
\multiput(56.75,76.5)(-.03373344371,-.03394039735){1208}{\line(0,-1){.03394039735}}
\put(16,35.5){\line(-1,0){.75}}
\multiput(56.75,76.5)(.03372374798,-.03392568659){1238}{\line(0,-1){.03392568659}}
\multiput(56.5,76.5)(-.0337136929,-.0840248963){482}{\line(0,-1){.0840248963}}
\multiput(56.25,76.25)(.0336842105,-.0826315789){475}{\line(0,-1){.0826315789}}
\put(16,35.25){\line(0,-1){31.5}}
\put(40.25,36.25){\line(0,-1){32}}
\put(72.25,37.25){\line(0,-1){32.25}}
\put(98.25,34.5){\line(0,-1){29.5}}
\put(15.791,34.541){\circle*{1.581}}
\put(40.541,35.541){\circle*{1.581}}
\put(72.791,36.541){\circle*{1.581}}
\put(98.291,33.791){\circle*{1.581}}
\put(57.291,75.291){\circle*{1.581}}
\put(16.041,4.041){\circle*{1.581}}
\put(40.541,4.291){\circle*{1.581}}
\put(72.541,5.291){\circle*{1.581}}
\put(98.541,5.291){\circle*{1.581}}
\put(48.059,35.809){\circle*{1.118}}
\put(54.559,36.059){\circle*{1.118}}
\put(61.059,36.059){\circle*{1.118}}
\put(47.309,4.559){\circle*{1.118}}
\put(53.809,4.809){\circle*{1.118}}
\put(60.309,4.809){\circle*{1.118}}
\put(57.25,82){\makebox(0,0)[cc]{$v_0$}}
\put(1.75,36.7){\makebox(0,0)[cc]{$v_1$}}
\put(31.5,37.75){\makebox(0,0)[cc]{$v_2$}}
\put(112,37.7){\makebox(0,0)[cc]{$v_{n-1}$}}
\put(85,38.75){\makebox(0,0)[cc]{$v_{n-2}$}}
\put(1.25,3){\makebox(0,0)[cc]{$u_1$}}
\put(33,3){\makebox(0,0)[cc]{$u_2$}}
\put(112,3){\makebox(0,0)[cc]{$u_{n-1}$}}
\put(85,3){\makebox(0,0)[cc]{$u_{n-2}$}}
\end{picture}
\end{center}

Then the minimum roman dominating distance matrix is
$$
  A_{RDd}(K^*_{1,n-1}) =
  \left[
\begin{array}{c|c|c}
2&\bf{1}_{n-1}&\bf{2}_{n-1}\\ \hline
\bf{1}'_{n-1}&2J_{n-1}-2I_{n-1} & 3J_{n-1}-2I_{n-1} \\ \hline
\bf{2}'_{n-1}&3J_{n-1}-2I_{n-1} & 4J_{n-1}-3I_{n-1}
\end{array}\right]_{2n-1\times 2n-1}$$
The characteristic equation of $A_{RDd}(K^{*}_{1,n-1})$ is $$(\rho^2+5\rho+2)^{n-2} [\rho^3-(6n-9)\rho^2-(n^2-7n+14)\rho+2n^2-3n-3]=0.$$
The sum of absolute values of the roots of $(\rho^2+5\rho+2)^2=0$ is $5n-10$ and the sum of absolute values of roots of the   $\rho^3-(6n-9)\rho^2-(n^2-7n+14)\rho+2n^2-3n-3=0$ is greater than or equal to $6n-9.$ Therefore, $$E_{RDd}(K^{*}_{1,n-1})\geq 11n-19.$$

By Cauchy’s bound for roots of  a polynomial, all the roots of $\rho^3-(6n-9)\rho^2-(n^2-7n+14)\rho+2n^2-3n-3=0$ lies in the closed interval $[-(M+1), M+1]$ where $M=2n^2-3n-3.$ Therefore, the sum of absolute values of these roots is bounded above by $3(2n^2-3n-2)=6n^2-9n-6.$
Thus, $$E_{RDd}(K^{*}_{1,n-1})\leq 5n-10+6n^2-9n-6=6n^2-4n-16.$$
\end{proof}

\begin{theorem}
For any  integer $k\geq 2,$  $E_{RDd}(S_k^0)=7k-6+\sqrt{k^2-4k+12}.$ 
\end{theorem}
\begin{proof}
Consider the crown graph $S_k^0$ with vertex set $V=X\cup Y$ where $X=\{v_1, v_2, \dots, v_k\}$ and $Y=\{w_1, w_2, \dots, w_k\}.$ The minimum roman dominating distance function is $f=(V_0, V_1, V_2)$ where $V_2=\{v_i, w_i\}$ for any $1\leq i\leq k,$ $V_1=\emptyset$ and $V_0=V\setminus V_2.$ Then the minimum roman dominating distance matrix is
$$A_{RDd}(S_k^0)=\left[
\begin{array}{c|c}
2J_k-D_i & J_k+2I_k \\ \hline
J_k+2I_k & 2J_k-D_i
\end{array}\right]_{2k\times 2k}$$
The Characteristic equation of $A_{RDd}(S_k^0)$ is
$$\rho^{2k-2}(\rho + 4)^{2k-2}[\rho^2 - (3k + 2)\rho + 6(k-1)] [(\rho^2 + (6 - k)\rho - 2k + 6] = 0.$$
Then the eigenvalues are $-4$ with multiplicity $k-2,$ 0 with multiplicity $k-2,$ $\frac{(-6+k)\pm \sqrt{k^2-4k+12}}{2}$ and $\frac{(3k+2)\pm \sqrt{9k^2-12k+28}}{2}.$
Hence the sum of absolute values of all eigenvalues is $7k-6+\sqrt{k^2-4k+12}.$ That is, $E_{RDd}(S_k^0)=7k-6+\sqrt{k^2-4k+12}.$
\end{proof}

\begin{theorem}
Let $G$ and $H$ be two disjoint graphs. Then $E_{RDd}(G\cup H)=E_{RDd}(G)+E_{RDd}(H).$
\end{theorem}
\begin{proof}
Let $A$ and $B$ be the minimum roman dominating distance matrix of $G$ and $H,$ respectively. Then the  minimum roman dominating distance matrix of $G\cup H$ is
$$A_{RDd}(G\cup H)=\left[
\begin{array}{c|c}
A &0\\ \hline
0 & B
\end{array}\right].$$
The characteristic polynomial of $A_{RDd}(G\cup H)$ is the product of characteristic polynomial of $A$ and $B.$
 Therefore, $E_{RDd}(G\cup H)=E_{RDd}(G)+E_{RDd}(H).$
\end{proof}
\section{ Properties of Eigenvalues of  Minimum Roman Dominating  Distance Matrix $A_{RDd}(G)$ }

In this section, we discussed the relation between the eigenvalues of the minimum roman dominating distance matrix $A_{RDD}(G)$ and the minimum roman dominating energy $\gamma_R$ of $G.$

\begin{theorem} 
Let $G=(V, E)$ be a graph and let $f=(V_0, V_1, V_2)$ be a $\gamma_R$-function. If $\rho_1, \rho_2, \dots, \rho_n$ are the eigenvalues of minimum roman dominating distance matrix $A_{RDd}(G),$ then
\begin{enumerate}
    \item[(i)] $\sum\limits_{i=1}^n \rho_i=\gamma_R(G)$
    \item[(ii)] $\sum\limits_{i=1}^n \rho_i^2=\gamma_R(G)+2m+2M$
 where $M =\sum\limits_{i<j,\, d(v_i, v_j)\neq 1} d(v_i, v_j)^2$ and $m = |E|.$
 \end{enumerate}
\end{theorem}
\begin{proof}
(i) The sum of the eigenvalues of $A_{RDd}(G)$ is the trace of $A_{RDd}(G).$ Therefore,
$$\sum\limits_{i=1}^n\rho_i=\sum\limits_{i=1}^n d(v_i, v_i)=2|V_2|+|V_1|=\gamma_R(G).$$

(ii) The sum of squares of the eigenvalues of $A_{RDd}(G)$ is trace of $(A_{RDd}(G))^2.$
Therefore,
\begin{align*}
    \sum_{i=1}^n \rho_i^2&=\sum_{i=1}^n\sum_{j=1}^n d(v_i, v_j)d(v_j, v_i)\\
    &=\sum\limits_{i=1}^n d(v_i, v_i)^2+\sum\limits_{i\neq j} d(v_i, v_j) d(v_j, v_i)\\
    &=\sum\limits_{i=1}^n d(v_i, v_i)^2+\sum\limits_{i<j} d(v_i, v_j)^2\\
    &= \gamma_R(G)+2\sum\limits_{i<j} d(v_i, v_j)^2\\
    &=\gamma_R(G)+2m+2M.
\end{align*}
\end{proof}
\begin{corollary}
Let $G$ be a  graph with diameter 2 and let $f=(V_0, V_1, V_2)$ be a $\gamma_R$-function. If $\rho_1, \rho_2,\dots,\rho_n$ are eigenvalues of minimum roman dominating distance matrix $A_{RDd}(G),$ then $$\sum\limits_{i=1}^n \rho_i^2=\gamma_R(G)+2(2n^2-2n-3m).$$
\end{corollary}
\begin{proof}
We know that in $A_{RDd}(G)$ there are $2m$ elements with 1 and $n(n - 1) - 2m$ elements with 2 and hence corollary follows from the above theorem.
\end{proof}

\section{Bounds for Minimum Roman Dominating Energy}
In this section, we discussed the bounds for minimum roman dominating energy. 

The proofs of the following Theorems are similar to the proofs in \cite{MRRaj}. 
\begin{theorem}
Let $G$ be a graph. If $f=(V_0, V_1, V_2)$ is a
$\gamma_R$-function  and $P = |det(A_{RDd}(G)|,$ then
$$\sqrt{(2m + 2M + \gamma_R) + n(n -1)P^{\frac{n}{2}}}
 \leq E_{RDd}(G) \leq \sqrt{n(2m + 2M + \gamma_R(G))}$$
where $\gamma_R$ is a roman domination number.
\end{theorem}

By Theorem\ref{a},we have the following Corollary
\begin{corollary}
Let $G$ be a  graph. If $f=(V_0, V_1, V_2)$ is a $\gamma_R$-function and $P = |det(A_{RDd}(G)|$ then
$$\sqrt{(2m + 2M + \gamma) + n(n -1)P^{\frac{n}{2}}}
 \leq E_{RDd}(G) \leq \sqrt{n(2m + 2M + 2\gamma(G))}$$
where $\gamma$ is a minimum domination number of $G.$
\end{corollary}

\begin{remark}
In \ref{b}, for the healthy spider graph $K^{*}_{1,n-1},$ $m=2n-2, M=(n-1)(19n-6)$ and $\gamma_{R}(K^{*}_{1,n-1})=n+1.$ Hence $\sqrt{n(2m + 2M + \gamma_R(G))}=\sqrt{(2n-1)(38n^2+31n-3)}>6n^2-4n-16.$
\end{remark}

\begin{theorem}
If $\rho_1(G)$ is the largest eigenvalue of a minimum roman dominating distance matrix
$A_{RDd}(G),$ then
$$\rho_1(G) \geq
\frac{2W(G) + \gamma_R(G)}{n}$$
where $W(G)$ is the Wiener index of $G.$
\end{theorem}
\begin{theorem}
Let $G$ be a graph of diameter 2 and $\rho_1(G)$ is the largest eigenvalue of a minimum roman
dominating distance  matrix $A_{RDd}(G),$ then
$$\rho_1(G) \geq
\frac{2n^2 - 2m - 2n + \gamma_R(G)}{n}.$$
\end{theorem}

\section*{Conclusion}
In this paper, we introduced the concept of minimum roman dominating  distance  energy $E_{RDd}(G)$ of a graph $G$ and computed minimum roman dominating  distance  energy of   complete, complete bipartite, crown, star  and healthy spider graphs. Also, we discussed the properties of eigenvalues of a minimum roman dominating distance  matrix $A_{RDd}(G).$ Finally, we derived the  Upper and lower bounds for $E_{RDd}(G).$

\end{document}